\documentclass[ 10 pt, peerreview]{ieeeconf}  

\IEEEoverridecommandlockouts                       
\usepackage{cite}
\usepackage{graphicx}

\usepackage{amsmath, amssymb, amsfonts} 
\usepackage{booktabs} 
\usepackage{mathtools}
\usepackage{dsfont}
\usepackage{commath}

\usepackage{microtype}
\usepackage{mathdots}
\usepackage{mathtools}
\usepackage{bm}

\usepackage{eucal}

\usepackage{tikz}
\usepackage{pgfplots}
\pgfplotsset{every axis/.append style={line width=1pt}}
\pgfplotsset{compat=newest}
\usetikzlibrary{shapes, fit, decorations.markings, matrix, plotmarks, positioning, fit, spy, patterns, shadows, calc, backgrounds, arrows.meta}
\pgfplotsset{every tick label/.append style={font=\tiny}}
\pgfplotsset{every axis/.append style={
                    label style={font=\footnotesize}
                    }}
\usepackage{algorithm} 
\usepackage{algpseudocode} 

\usepackage[font=footnotesize,caption=false]{subfig}
\usepackage{graphicx}
\usepackage[abs]{overpic}

\usepackage[nolist]{acronym}
\usepackage{siunitx}
\usepackage{booktabs}
\usepackage{lipsum}

\usepackage{enumerate}
\newcounter{subeqn} %

\newtheorem{theorem}{Theorem}
\newtheorem{corollary}{Corollary}
\newtheorem{lemma}[theorem]{Lemma}
\newtheorem{proposition}{Proposition}
\newtheorem{assumption}{Assumption}
\newtheorem{definition}{Definition}

\newcommand{\ts}[1]{{\textnormal{#1}}}
\newcommand{\ie}{\emph{i.e.,}\ }

\newcommand{\mc}{\mathcal}
\newcommand{\mb}{\mathbf}

\newcommand{\mbb}{\mathbb}

\newcommand{\Nset}{\mathbb{N}}
\newcommand{\Rset}{\mathbb{R}}

\newcommand{\mf}{\mathfrak}
\newcommand{\mn}{\mathrm}

\newcommand{\relu}[1]{\phi(#1)}

\newcommand{\Pre}{\operatorname{Pre}}
\newcommand{\im}{\operatorname{im}}
\newcommand{\rank}{\operatorname{rank}}
\newcommand{\ini}{\mathrm{ini}}
\newcommand{\eqm}{\mathrm{eq}}

\newif\ifproves
\provestrue

\allowdisplaybreaks

\begin{acronym}
  \acro{DeePC}{Data Enabled Predictive Control}
  \acro{MPC}{Model Predictive Control}
  \acro{PE}{Persistently Exciting}	
  \acro{LMI}{Linear Matrix Inequalities}
  \acro{NN}{Neural Network}
  \acro{MPI}{Maximal Positive Invariant}
  \acro{MCI}{Maximal Control Invariant}
  \acro{PWA}{Piece-wise Affine}
\end{acronym}
\usepackage{cite}
\title{Feasibility certificates for DeePC: A set invariance approach}
\author{Daniel O'Neill$^{1}$, Ross Drummond,$^{1}$ and Pablo R. Baldivieso-Monasterios$^{1}$\thanks{$^{1}$ Daniel O'Neill, Ross Drummond, and Pablo R. Baldivieso-Monasterios are with the School of Electric and Electronic Engineering, University of Sheffield, Sheffield, UK. Emails: {\tt\footnotesize \{doneil6,ross.drummond, p.baldivieso\}@sheffield.ac.uk}.}
}
\definecolor{col_1}{rgb}{0.19054, 0.01184, 0.68504}
\definecolor{col_2}{rgb}{0.91433, 0.55532, 0.12457}
\definecolor{col_3}{rgb}{0.37363, 0.14399, 0.00538}
\definecolor{col_4}{rgb}{0.57815, 0.63390, 0.73813}
\definecolor{col_5}{rgb}{0.11766, 0.42406, 0.32378}
\definecolor{col_6}{rgb}{0.74572, 0.59932, 0.52154}
\definecolor{col_7}{rgb}{0.39925, 0.32996, 0.46159}
\definecolor{col_8}{rgb}{0.69631, 0.42036, 0.18407}
\definecolor{col_9}{rgb}{0.29856, 0.16597, 0.35122}
\definecolor{col_10}{rgb}{0.14331, 0.66453, 0.98936}
\definecolor{col_11}{rgb}{0.16172, 0.24576, 0.21297}
\definecolor{col_12}{rgb}{0.09245, 0.51449, 0.71787}
\definecolor{col_13}{rgb}{0.39726, 0.12324, 0.98385}
\definecolor{col_14}{rgb}{0.53408, 0.48188, 0.50616}
\definecolor{col_15}{rgb}{0.53623, 0.28082, 0.09604}
\definecolor{i_col_1}{rgb}{0.50634, 0.68504, 0.01184}
\definecolor{i_col_2}{rgb}{0.12457, 0.48358, 0.91433}
\definecolor{i_col_3}{rgb}{0.00538, 0.23502, 0.37363}
\definecolor{i_col_4}{rgb}{0.73813, 0.68238, 0.57815}
\definecolor{i_col_5}{rgb}{0.42406, 0.11766, 0.21794}
\definecolor{i_col_6}{rgb}{0.52154, 0.66793, 0.74572}
\definecolor{i_col_7}{rgb}{0.39230, 0.46159, 0.32996}
\definecolor{i_col_8}{rgb}{0.18407, 0.46001, 0.69631}
\definecolor{i_col_9}{rgb}{0.21863, 0.35122, 0.16597}
\definecolor{i_col_10}{rgb}{0.98936, 0.46815, 0.14331}
\definecolor{i_col_11}{rgb}{0.24576, 0.16172, 0.19451}
\definecolor{i_col_12}{rgb}{0.71787, 0.29584, 0.09245}
\definecolor{i_col_13}{rgb}{0.70983, 0.98385, 0.12324}
\definecolor{i_col_14}{rgb}{0.48188, 0.53408, 0.50980}
\definecolor{i_col_15}{rgb}{0.09604, 0.35144, 0.53623}

\begin{document}
\bstctlcite{nourl}
\bstctlcite{nodoi}
\maketitle

\thispagestyle{empty} 
\pagestyle{empty}

\begin{abstract}      
  \ac{DeePC} is posed on a window of past inputs and outputs, so its feasible set is a set of histories and is not available in closed form. We write the \ac{DeePC} solution as an implicit ReLU network whose weights are built from Hankel matrices and whose activation patterns are the critical regions of the underlying quadratic program. This gives, directly from the Hankel matrices, the feasible set of DeePC as a polyhedral projection, its explicit piecewise-affine control law, and an exact linear model of the history window. We compute three invariant sets of histories: the maximal positively invariant set of the unsaturated linear branch, the region of recursive feasibility of DeePC without terminal ingredients, and the maximal control invariant set of admissible histories. We show these sets are nested and provide numerical examples.  
\end{abstract}


\section{Introduction }
\label{sec:introduction}
Consider an unknown discrete linear time invariant system $\Sigma = (\mc{U},\mc{Y},\mc{B})$ relating an input space $\mc{U}$ to an output space $\mc{Y}$ by $\mc{B}\subset\mc{U}\times\mc{Y}$ and a data set $\mc{S}$ gathered from the system, \ie  $\mc{S} = \{(u_i^{d},y_i^d)\}_i^{T}$ is a finite sequence of input/output pairs. The celebrated Willems' Lemma \cite{Willems2005} allows us to describe the linear system $\Sigma$ using only information $\mc{S}$, provided the input sequence is \ac{PE}, which will be defined below.
\begin{definition}[see \cite{DePersis2020}]
  For $T>0$ and  $v = \{v_i\}_i^{T}\subset\Rset^{m}$,
  \begin{enumerate}[i)]
  \item the depth-$L$ block-Hankel matrix is
\begin{align*}
    \mc{H}_L(v) := \begin{bmatrix}
v_1 & v_2 & \cdots & v_{T-L+1} \\
v_2 & v_3 & \cdots & v_{T-L+2} \\
\vdots & \vdots & \ddots & \vdots \\
v_L & v_{L+1} & \cdots & v_{T}
\end{bmatrix},
\end{align*}
\item the sequence $v\subset\Rset^m$ is \emph{Persistently Exciting} of order $L$ if $\ts{rank}(\mc{H}_L(v)) = mL$.
  \end{enumerate}
  \label{def:hank_PE}
\end{definition}

Definition~\ref{def:hank_PE}$-ii)$ together with \cite[Lemma 2]{DePersis2020}, which is a restatement of \cite[Thm 1]{Willems2005}, states that for a data set $\mc{S}$ such that $u^d$ is \ac{PE} of order $n+t$, then any $L-$long sequence pair $(u,y) =\{(u_0,y_0),\ldots,(u_{L-1},y_{L-1})\} \in \mc{B}$ can be written as a linear combination of the elements of $\mc{S}$, \ie for some $g\in\Rset^{L+1}$
\begin{equation}
  \label{eq:Hankel_basis}
    \begin{bmatrix}
    u\\y
  \end{bmatrix} =
  \begin{bmatrix}
    \mc{H}_L(u^d)\\\mc{H}_L(y^d)
  \end{bmatrix}g.
\end{equation}

Equation \eqref{eq:Hankel_basis} is the cornerstone of data-driven controllers. The Hankel matrix of the data set $\mc{S}$ generates all possible trajectories. The relation between $\mc{H}_L$ and the trajectories $\{u,y\}$ was established by contradiction, obscuring the relationship between the linear system $\Sigma$ and $\mc{S}$. In \cite{Markovsky2023}, the authors overcome this conceptual barrier by providing constructive proofs for the results in \cite{Willems2005}. Furthermore, \cite{Lopez2022} extends the framework to the continuous-time case.
For more in-depth reviews of some recent developments in data-driven predictive control, we refer the reader to the excellent reviews \cite{Markovsky2021, Berberich2025}.  

\subsection{Data-enabled predictive control}
\label{sec:DD_QP}
\ac{DeePC}, initially posed by \cite{Coulson2019} and rigorously analysed in \cite{Berberich2021c}, exploits the relationship~\eqref{eq:Hankel_basis} and includes it within an optimisation problem. However, \ac{DeePC}, unlike traditional \ac{MPC}, is inherently an input/output controller. Since the current state is unavailable, the current operating condition is represented by a finite history window of the $T_{\mathrm{ini}}>0$ most recent measured inputs and outputs, $\xi = (u_{\mathrm{ini}}, y_{\mathrm{ini}})\in \Rset^{(m+p)T_{\mathrm{ini}}}$, where $u_\mn{ini} = (u_{k-T_\mn{ini}},\ldots,u_{k-1})\subset\Rset^{mT_\mn{ini}}$ and $y_\mn{ini} = (y_{k-T_\mn{ini}},\ldots,y_{k-1})\subset\Rset^{pT_\mn{ini}}$. Furthermore, $T_\mn{ini} >0$ satisfies the following assumption:
\begin{assumption}[Sufficient initial window length] 
  The initial trajectory length satisfies $T_{\mathrm{ini}} \geq \ell,$ where $\ell$ is the system observability index\footnote{For a system $x^+ = Ax + Bu$ with output $y = Cx$ and $x\in\Rset^n$, the observability index $\ell\in\Nset$ is the minimum number such that the observability matrix $[C^\top~~(CA)^\top~~\cdots~~(CA^{\ell-1})^\top]^\top$ has column rank $n$.}.
  \label{assump:Tini}
\end{assumption}

Under Assumption \ref{assump:Tini}, $\xi$ contains the information about the system's internal condition required for prediction. The $N-$horizon \ac{DeePC} optimises over future input output predictions, $ \mb{u} = (u_0 ,\ldots, u_{N-1} )$ and $ \mb{y} = (y_0 ,\ldots, y_{N-1} )$ respectively, to regulate the system $\Sigma$ towards a desired reference $ \mb{r} = (r_0 ,\ldots, r_{N-1})$. The associated quadratic program is 
\begin{subequations}
  \begin{align}
\mbb{P}(\mb{r},\xi)\colon & \min_{g, \mb{u}, \mb{y}} J(\mb{y},\mb{u},\mb{r}) \\
    \text{s.t.} \quad &  u_{k} \in \mbb{U}, \quad k = 0, \dots, N-1, \label{eq:DeePC_u}\\
                          & y_{k} \in\mbb{Y}, \quad k = 0, \dots, N-1,\label{eq:DeePC_y}\\
                          & \begin{bmatrix} \mc{U}_p \\ \mc{Y}_p\end{bmatrix} g
  = \xi,\,\,\begin{bmatrix} \mc{U}_f \\ \mc{Y}_f \end{bmatrix} g
                            = \begin{bmatrix} \mb{u} \\ \mb{y} \end{bmatrix},\label{eq:fund-lemma}
  \end{align}
  \label{eq:DeePC_QP}
\end{subequations}
where 
\begin{align*}
\begin{bmatrix} \mc{U}_p\\ \mc{U}_f \end{bmatrix} = \mc{H}_L(u^d) ,\qquad
\begin{bmatrix} \mc{Y}_p\\ \mc{Y}_f \end{bmatrix} = \mc{H}_L(y^d) .
\end{align*}

The relations in \eqref{eq:fund-lemma} ensure consistency with the measured input/output history and future input/output predictions. Consequently, the measured history $\xi$ serves as the known initial condition for the behavioural predictor without explicitly reconstructing the plant state. The cost, in the standard framework, penalises deviations from a reference: $ J (\mb{y},\mb{u},\mb{r}) = \sum_{k=0}^{N-1} \left( \|y_{\mathrm{fut}, k} - r_k\|_Q^2 + \|u_{\mathrm{fut}, k} - u_{\mathrm{ref}, k}\|_R^2 \right)$. In \cite{Dorfler2023}, the authors modified problem \eqref{eq:DeePC_QP} by adding regularisation terms to the cost to ensure the optimisation problem is well-posed. Slack variables in \eqref{eq:fund-lemma}, common in practice \cite{Coulson2019}, soften the consistency constraint, so that \eqref{eq:DeePC_QP} becomes feasible for histories that no trajectory of $\Sigma$ generates; they are excluded in this paper. The authors of \cite{Wang2026} and \cite{Berberich2020a} reframe the \ac{DeePC} problem into a min-max robust formulation that optimises worst-case performance over an uncertainty set derived from the non-unique solutions to the behavioural data equations. The constraints of \eqref{eq:DeePC_QP} satisfy the following assumption.
\begin{assumption}
  The sets $\mbb{U}\subset\Rset^m$ and $\mbb{Y}\subset\Rset^p$ are convex and compact. In addition, $\mbb{U}$ contains the origin in its non-empty interior.
  \label{assump:constraints}
\end{assumption}

Stability of the origin is ensured via the addition of suitable terminal ingredients, similar to model-based \ac{MPC}. \cite{Berberich2021a} includes these terminal sets in the \ac{DeePC} formulation; \cite{Berberich2022} extends these ideas to a nonlinear setting, ensuring tracking of a given target. In contrast, \cite{Faro2025} modifies~\eqref {eq:DeePC_QP} to allow updates to the Hankel matrices by imposing additional constraints on $\mb{u}$. We refer the reader to \cite{Verheijen2023} for an in-depth review of \ac{DeePC}. 

\subsection{Set Invariance and Feasibility in \ac{DeePC}: Current Limitations and Contributions}

Recursive feasibility is a key issue for predictive-based controllers because it ensures that solutions to \eqref{eq:DeePC_QP} exist along the closed-loop trajectories. In nominal linear time-invariant settings with \ac{PE} data, direct \ac{DeePC} yields an initial feasible set that is theoretically equivalent to that of standard model-based output-feedback MPC \cite{Coulson2019,Mattsson2024}. A similar result was obtained by \cite{Watson2025}, who studied feasibility in the case of partial information (a hybrid between data-driven and model-based approaches). To enforce recursive feasibility and closed-loop stability, \cite{Berberich2021,Berberich2025} introduced terminal equality constraints and data-driven terminal invariant sets, which guarantee a valid shifted trajectory at each time step while expanding the controllable region. Under measurement noise and disturbances, tube-based and robust min-max formulations preserve recursive feasibility by tightening state/input constraints and introducing initial-condition slack variables \cite{Huang2023}. Lastly, \cite{Kladtke2024a} proposes a generalisation of finite control sets \ac{MPC} to \ac{DeePC}.

Computing feasible sets in data-driven predictive control requires computing a set for histories $\xi$ rather than an initial state. This poses an issue when projecting dataset-dependent trajectory variables. This input-state framework has been used to certify linear control laws by constructing data-driven invariant sets \cite{Bisoffi2023,Bisoffi2020} or control barrier functions directly from noise-corrupted trajectories \cite{Bajelani2024a, Liao2026}. Furthermore, maintaining recursive feasibility without shrinking the admissible feasible domain under noisy conditions requires specialised set-expansion techniques and robust safety-critical relaxations \cite{Bajelani2025, Liao2026}. In this paper, we will focus on two fronts: \ac{DeePC} and set-theoretic computations for data-driven predictive controllers. We propose a method for computing the feasible set of \eqref{eq:DeePC_QP} that avoids projecting out the dataset-dependent variable $g\in\Rset^{T-L+1}$: $g$ is eliminated in closed form, and only the $mN$ future inputs remain. We exploit the structure of \eqref{eq:DeePC_QP} to formulate an implicit \ac{NN} \cite{ElGhaoui2021,Drummond2024} whose weights are computed from the Hankel matrices and whose activation patterns are the critical regions of \eqref{eq:DeePC_QP}. 
The contributions of this paper are:
\begin{enumerate}[i)]
\item An interpretation of the \ac{DeePC} problem in terms of implicit \ac{NN}s whose activation patterns are the critical regions of \eqref{eq:DeePC_QP}. This representation allows us to characterise the region where constraints are not triggered. Our first invariant set, akin to the model-based \ac{MPI} set and the data-driven version from \cite{Bajelani2024a}, is in Proposition~\ref{prop:dd-oinf}.  
\item A characterisation of the \ac{DeePC} feasible region in the augmented state $\xi$, of the explicit \ac{PWA} \ac{DeePC} law, and of an exact linear model of the history window, all computed from data. Proposition~\ref{prop:crit} partitions the history space into regions where a linear gain can be applied and shows that these constitute the feasible set of \ac{DeePC}.
\item Construction of data-driven invariant sets: the invariant set of the unsaturated branch, the region of recursive feasibility of \ac{DeePC} without terminal ingredients, and the maximal control invariant set of admissible histories, together with their ordering. Theorem~\ref{thm:hierarchy} and its corollary show that these invariant sets nest.
\end{enumerate}
These results complement existing feasibility certificates in the \ac{DeePC} literature and relate them to their model-based counterparts.

\textit{\underline{Notation}:} The $n\times n$ identity matrix is $\mbb{I}_n$. For a matrix $A\in\Rset^{n\times m}$, $A^{\dagger}$ is the Moore--Penrose pseudo inverse, $\im A$ its range, and, for an index set $\mc{A}$, $A_\mc{A}$ collects the rows indexed by $\mc{A}$, $A_{\mc{A}\mc{B}}$ the corresponding sub matrix, and $\bar{\mc{A}}$ is the complement of $\mc{A}$. $\rho(\cdot)$ is the spectral radius and $\relu{v} = \max(v,0)$, taken component-wise, is the ReLU. Inequalities between vectors hold component-wise.

\begin{definition}[Implicit neural network]
  An implicit \ac{NN} $\mc{N}\colon\Rset^{n_z}\to\Rset^{n_o}$ is defined by the solution of the following system of equations
 \begin{subequations}\label{eqn:mpc-to-nn}
      \begin{align}
        \zeta= & W\phi(\zeta) + Yz + b,\\
        \mc{N}(z) = &W_f\phi(\zeta) + Y_f z + b_f,
      \end{align}
    \end{subequations}
    for an activation $\relu{\cdot}$, weights $W\in\Rset^{n_c\times n_c}$, $Y\in\Rset^{n_c\times n_z}$, $W_f\in\Rset^{n_o\times n_c}$, $Y_f\in\Rset^{n_o\times n_z}$, and biases $b\in\Rset^{n_y}$, $b_f\in\Rset^{n_o}$.
    \label{def:Im_NN}
\end{definition}

\begin{definition}[\ac{MPI}]
  For $f\colon\mc{D}\to\mbb{V}$ and $\mc{S}\subseteq\mc{D}$, the \ac{MPI} subset of $\mc{S}$ is $\mc{O}(f,\mc{S}) := \{\xi\in\mc{S}\colon f^k(\xi)\in\mc{S},~\forall k\in\Nset\}$ where $f^k(\cdot)$ is the repeated application of the map $f(\cdot)$.
  \label{def:MPI}
\end{definition}

\section{Implicit neural network representation}
\label{sec:impl-neur-netw}

In this section, we express~\eqref{eq:DeePC_QP} using an implicit neural network developed in \cite{Drummond2024}. We condense \eqref{eq:DeePC_QP} using the following lemma:
\begin{lemma}[Condensed \ac{DeePC} QP]
Suppose Assumptions \ref{assump:Tini} and \ref{assump:constraints} hold. Suppose, in addition, $Q\succeq 0$, $R\succ 0$, and given $\mc{N}$ the basis of $\ts{ker}(A_{\mn{eq}})$ where $ A_{\mn{eq}} = [ \mc{U}_p^\top \,\,\, \mc{Y}_p^\top ]^\top$, $\mc{U}_f\mc{N}$ has full row rank. The QP~\eqref{eq:DeePC_QP} is equivalent to
  \begin{subequations}
    \begin{align}
      \min\,\,& \tfrac{1}{2}\mb{u}^\top H\mb{u} + z^\top F^\top\mb{u},\label{eq:cond_cost}\\
    \ts{s.t.}\,\,\,&G\mb{u} \leq Sz + w,\label{eq:cond_cons}
    \end{align}
    \label{eq:DeePC_QP_cond}
  \end{subequations}
  where $z = (\xi,\mb{r})$, $H\in\Rset^{mN\times mN}$ such that $H\succ 0$, $F \in\Rset^{(m+p)T_\mn{ini}+Np\times Nm}$, $G\in\Rset^{(n_u+n_y)N\times Nm}$, $S\in\Rset^{(n_u+n_y)N\times (m+p)T_\mn{ini}+Np}$, and $w\in\Rset^{(n_u+n_y)N}$.
  \label{lem:QP_cond}
\end{lemma}
\begin{proof}
  First, we eliminate the behavioural coefficient $g$; in \eqref{eq:fund-lemma}, its role splits into past and future components. The future trajectory is $\mb{u} = \mc{U}_fg$ and $\mb{y} = \mc{Y}_f g$. Replacing this in the cost yields\[
  \begin{split}
    J(g) = g^\top(\mc{U}_f^\top(\mbb{I}_N\otimes R)&\mc{U}_f + \mc{Y}_f^\top(I_N\otimes Q)\mc{Y}_f)g + \\
    -& 2\mb{r}^\top(I_N\otimes Q)\mc{Y}_fg + \mb{r}^\top(I_N\otimes Q)\mb{r}.  
  \end{split}\]
The consistency conditions are equality constraints $A_{\mn{eq}}g=\xi $. Solving these yield $g =A_{\mn{eq}}^{+}\xi + \mc{N}c$ with $c\in\Rset^{\ts{dim}(\ts{ker}(A_\mn{eq}))}$. By assumption, $\mc{U}_f\mc{N}$ has full row rank, then $V=\mathcal N(\mathcal U_f\mathcal N)^\dagger$, yields $\mathcal U_fV=\mbb{I}_{mN}$. From $\mb{u}=\mc{U}_fg$, we obtain $\mb{u} = \mc{U}_fA_\mn{eq}^\dagger \xi + \mc{U}_f\mc{N}c$ which leads to $c = (\mc{U}_f\mc{N})^\dagger(\mb{u} - \mc{U}_fA_\mn{eq}^\dagger\xi)$. The behavioural coefficient vector $g$ can be expressed as
\begin{equation}
g = \underbrace{(\mbb{I}_{T}-V\mc{U}_f)A_{\mn{eq}}^\dagger}_{A_0}\xi + V\mb{u} + \mc{N}c_0, \qquad c_0\in\ker\mc{U}_f\mc{N}.
\label{eq:g-parameterisation} 
\end{equation}
The matrix $V$ satisfies $\mc{U}_fA_0 = 0$. Since $\mc{N}c_0\in\ker[A_\mn{eq}^\top~\mc{U}_f^\top]^\top\subseteq\ker\mc{Y}_f$ (zero input generates zero outputs), the prediction $\mb{y} = \mc{Y}_f(A_0\xi+V\mb{u})$ does not depend on $c_0$; as $J$ and \eqref{eq:DeePC_u}--\eqref{eq:DeePC_y} depend on $g$ only through $(\mb{u},\mb{y})$, the QP \eqref{eq:DeePC_QP} becomes a problem in $\mb{u}$ alone.
From Assumption~\ref{assump:constraints}, we have that there exist $P_y \in\Rset^{n_y\times p}\,q_y\in\Rset^{n_y}$ and $P_u \in\Rset^{n_u\times m}\,q_u\in\Rset^{n_u}$ with $n_y>0$ and $n_u>0$ that approximate the sets $\mbb{Y}$ and $\mbb{U}$. Substituting $\mathbf u=\mathcal U_fg$ and $\mathbf y=\mathcal Y_fg$ into the polytopic constraints and stacking the constraints over the prediction horizon gives the prediction constraints
$G\theta \leq w + [S_\xi\,\,S_r]z$ where $G = [(\mbb{I}_N\otimes P_u)^\top\,\,V^\top\mc{Y}_f^\top(\mbb{I}_N\otimes P_y)^\top]^\top$, $w = [(\mb{1}_N\otimes q_u)^\top\,\,(\mb{1}_N\otimes q_y)^\top]^\top$, and
\[S = [S_\xi\,\, S_r] = \begin{bmatrix} \mb{0} & \mb{0}\\ -(\mbb{I}_N\otimes P_y)\mc{Y}_fA_0 & \mb{0} \end{bmatrix}.\] Lastly, the cost matrices are the following:\[
  \begin{split}
    H =& V^\top\underbrace{(\mc{U}_f^\top(\mbb{I}_N\otimes R)\mc{U}_f + \mc{Y}_f^\top(I_N\otimes Q)\mc{Y}_f)}_{\tilde{H}}V,\\
    F = &\begin{bmatrix}
            F_\xi \\ F_r
          \end{bmatrix} = \begin{bmatrix}
            A_0^\top\tilde{H}\\ (I_N\otimes Q)\mc{Y}_f
          \end{bmatrix}V.
  \end{split}
\]
Finally, $\mc{U}_fV = \mbb{I}_{mN}$ gives $H = \mbb{I}_N\otimes R + V^\top\mc{Y}_f^\top(\mbb{I}_N\otimes Q)\mc{Y}_fV\succeq \mbb{I}_N\otimes R\succ 0$.
\end{proof}

The condition of $\mc{U}_f\mc{N}\in\Rset^{mN\times \ts{dim}\ker A_\mn{eq}}$ being full rank imposes a restriction on the data sequence length, \ie $T_\mn{ini}(m+p) + mN \leq T$ following the dimension theorem. With the \ac{DeePC} QP in condensed form, we can express its solution using the implicit \ac{NN}s from Definition~\ref{def:Im_NN}. This can be summarised in the following lemma, which has been adapted to our setting from \cite[Theorem 6]{Drummond2024}
\begin{lemma}[Adapted from \cite{Drummond2024}]\label{lem:mpc-to-nn}
  Suppose Lemma~\ref{lem:QP_cond} holds. The \ac{DeePC} problem mapping a history $z = (\xi,\mb{r})$ to inputs $\mb{u}$ written in the compact form
\begin{subequations}\label{eq:QP1}\begin{align}
  \mb{u}^*(z) = \arg \min &~\tfrac{1}{2} \mb{u}^\top H \mb{u} + z^\top  F \mb{u},\\
  \ts{s.t.~} & G \mb{u} \leq S z  + w,\label{eq:QP_com}
\end{align}
\end{subequations}
where $H \succ 0$ is a positive definite Hessian, and $F \in \mathbb{R}^{n \times m\,N}$. This mapping can be expressed as an implicit neural network
\begin{subequations}\label{eq:NN_controller}\begin{align}
    \zeta &= W \relu{\zeta}+Yz+ b,\label{eq:NN_implicit}\\
    \mb{u}^*(z) &= W_f\relu{\zeta} +Y_fz,\label{eq:NN_law}      
\end{align}\end{subequations}
with $\relu{\cdot}$ being the ReLU activation and with analytic expressions for the weights and biases given by
\begin{subequations}
  \begin{align}
    W &= \mbb{I}_{n_y} - GH^{-1}G^\top,\, W_{f} = -H^{-1}G^\top,\quad b = -w , \\
    Y &= -S-GH^{-1}F^\top,\, Y_{f} = - H^{-1}F^{\top}, \quad b_{f} = 0.
  \end{align}\label{eq:Wandb}
\end{subequations}
\end{lemma}

\section{Zero-Activation Region and Feedback Law}
\label{sec:dd-unsaturated}
In this section, we analyse the closed-loop behaviour inside the linear branch. The zero-activation condition provides the simplest affine branch of the constrained controller \eqref{eq:NN_law}. When $\relu{\zeta} = 0$, the implicit \ac{NN} output is
\begin{equation}
\mb{u}^*(\xi,\mb{r}) = -H^{-1}(F_\xi\xi + F_r\mb{r}).
\label{eq:theta-unsat}
\end{equation}
The associated behavioural parameter is \(g^* = \Gamma\xi - VH^{-1}F_r\mb{r}\) with $\Gamma = (A_0- VH^{-1}F_\xi)$. Replacing this control law into the constraints, we have the admissible set $\mc{X}$ for initial histories:
\[\begin{split}
  \begin{bmatrix} (\mbb{I}_N\otimes P_u)\mc{U}_\mn{ini}\\(\mbb{I}_N\otimes P_y)\mc{Y}_\mn{ini} \end{bmatrix} \Gamma\xi \leq &\begin{bmatrix} \mb{1}_N\otimes q_u\\ \mb{1}_N\otimes q_y \end{bmatrix} \\  
  - & \begin{bmatrix} (\mbb{I}_N\otimes P_u)\mc{U}_f\\(\mbb{I}_N\otimes P_y)\mc{Y}_f \end{bmatrix}VH^{-1}F_r\mb{r}.
\end{split}
\]
Furthermore, the zero-activation branch is polyhedral.
\begin{proposition}[Zero-activation admissible region]
  Suppose Assumptions~\ref{assump:Tini} and \ref{assump:constraints} hold. Fix $\mb{r} \in\Rset^{pN}$. On the set
  \begin{equation}
    \Xi_0(\mb{r}) = \left\{ \xi\colon  S_\xi\xi + w \geq - S_r\mb{r} \right\},
    \label{eq:Xi0}
  \end{equation}
  the solution of \eqref{eq:NN_implicit} at any $\xi\in\Xi_0(\mb{r})\cap\mc{X}$ satisfies $\relu{\zeta} = 0$.
  \label{prop:ZAAR}
\end{proposition}
\begin{proof}
  The solution of the nonlinear equation~\eqref{eq:NN_implicit} can be found using fixed point iterations, \ie $\zeta_{k+1} = W\zeta_k + Yz + b$ where $z = (\xi,\mb{r})$. From Lemma~\ref{lem:mpc-to-nn}, if $\xi \in \Xi_0(\mb{r})$, then $Yx + b \leq 0$. Consider, an initial $\zeta_0 \in\Rset^{n_c}$ such that w.l.o.g $\zeta_0 \geq 0$. Then, $W\relu{\zeta_0} = W{\zeta_0}$ and $Yz + b \in\Rset^{n_c}_{\leq 0}$ lead to an iterate $\zeta_1$ approaching the nonpositive orthant. As long as $\zeta_k\not\in\Rset^{n_c}_{\leq 0}$, then $\relu{\zeta_k} \neq 0$. However, since $Yz+b$ is constant and negative, the iterates converge to $\Rset^{n_c}_{\leq 0}$, \ie $\relu{\zeta_k}\to 0$.
\end{proof}

The set $\Xi_0(\mb{r})$ is not necessarily bounded. Proposition~\ref{prop:dd-oinf} below shows that the rows of $\mc{X}$ restrict only the initial history. 

From now on we consider regulation, $\mb{r} = 0$, so that $z = \xi$, $F = F_\xi$ and $S = S_\xi$. We write $\Xi_0 = \Xi_0(\bm{0})$, 
\begin{equation}
  \hat S := S_\xi + GH^{-1}F_\xi^\top,\qquad M := GH^{-1}G^\top,
  \label{eq:Shat}
\end{equation}
so that $Y = -\hat S$, and $E_1 := [\mbb{I}_m~~0]\in\Rset^{m\times mN}$ selects the input $u$ applied by the receding-horizon controller.

\subsection{History Window Closed-Loop Dynamics}
\label{sec:dd-unsaturated-cl}
The closed-loop analysis requires describing the system evolution in terms of the current $\xi$ and $u$. Let $\Pi$ be the shift that removes the oldest input/output pair from $\xi$ and $E$ the matrix that inserts a new pair, so that $\xi^+ = \Pi\xi +E(u_k,y_k)$, and let $\mc{U}_f^{(1)}$, $\mc{Y}_f^{(1)}$ be the first $m$ and $p$ rows of $\mc{U}_f$, $\mc{Y}_f$.
\begin{lemma}[History dynamics]
  Suppose Assumption~\ref{assump:Tini} hold, let $\mbb{V} = \im{A_\eqm}$, and  let $[L_\xi~~L_u] = \mc{Y}_f^{(1)}\Theta^\dagger$ with $\Theta = [A_\eqm^\top~~\mc{U}_f^{(1)\top}]^\top$. Every trajectory of $\Sigma$ satisfies $y = L_\xi\xi+L_uu$ and
  \begin{equation}
    \xi^+ = A_\xi\xi + B_\xi u,\quad A_\xi := \Pi + E\begin{bmatrix}0\\L_\xi\end{bmatrix},\quad B_\xi := E\begin{bmatrix}\mbb{I}_m\\L_u\end{bmatrix}.
    \label{eq:history-model}
  \end{equation}
  Moreover, $A_\xi\mbb{V}\subseteq\mbb{V}$ and $\im B_\xi\subseteq\mbb{V}$.
  \label{lem:history}
\end{lemma}
\begin{proof}
For $\bar g\in\ker\Theta$, by Assumption~\ref{assump:Tini}, we have that $\mc{Y}_f^{(1)}\bar g = 0$. Hence $\mc{Y}_f^{(1)} = \mc{Y}_f^{(1)}\Theta^\dagger\Theta$. By \eqref{eq:Hankel_basis}, every trajectory segment satisfies $(\xi,u,y) = (A_\eqm g, \mc{U}_f^{(1)}g,\mc{Y}_f^{(1)}g)$ for some $g$, so $y = \mc{Y}_f^{(1)}\Theta^\dagger\Theta g = L_\xi\xi+L_uu$, and \eqref{eq:history-model} follows. For $\xi\in\mbb{V}$ and any $u$, $\xi^+$ is a history of $\Sigma$, hence in $\mbb{V}$; taking $u = 0$ gives $A_\xi\mbb{V}\subseteq\mbb{V}$ and then $\im B_\xi\subseteq\mbb{V}$.
\end{proof}

On $\Xi_0$, the closed-loop is linear
\begin{equation}
  \xi^+ = A_{\mathrm{cl}}\xi,\qquad A_{\mathrm{cl}} := A_\xi - B_\xi H^{-1}F_\xi,
  \label{eq:Acl}
\end{equation}
built from the Hankel matrices only. On $\mbb{V}$, $A_{\mathrm{cl}}$ coincides with the map $\Pi + E[\mc{U}_f^{(1)\top}~\mc{Y}_f^{(1)\top}]^\top(A_0-VH^{-1}F_\xi)$ obtained from the first predicted input/output pair of the QP \eqref{eq:DeePC_QP}.

Membership of $\Xi_0$ only guarantees that the unconstrained prediction is admissible at the current sampling instant. The \ac{MPI} $\mc{O}_\infty = \mc{O}(A_\mn{cl},\Xi_0)$ from Definition~\ref{def:MPI} is the set that can certify the closed-loop dynamics inside $\Xi_0$. The computation of $\mc{O}_\infty$ follows the algorithm of Gilbert and Tan \cite{Gilbert1991}: $\mc{P}_0 := \Xi_0$ and
\begin{equation}
  \mc{P}_{k+1} := \mc{P}_k\cap\{\xi\colon \hat S A_{\mathrm{cl}}^{k+1}\xi+w\geq0,~A_{\mathrm{cl}}^{k+1}\xi\in\mc{X}\}.
  \label{eq:gt-dd}
\end{equation}
This process and its properties are summarised in the following result:
\begin{proposition}[Unsaturated invariant set]
  Suppose Assumptions~\ref{assump:Tini}--\ref{assump:constraints} hold, $\mb{r}=0$, and $A_{\mn{cl}}$ restricted to $\mbb{V}$ is Schur. Then:
  \begin{enumerate}[i)]
  \item there is a finite $k^\star$ with $\mc{P}_{k^\star+1}=\mc{P}_{k^\star} = \mc{O}_\infty$, and $\mc{O}_\infty$ is a polytope;
  \item for $\xi\in\mc{O}_\infty$, The QP \eqref{eq:DeePC_QP} is feasible along the closed loop trajectory of \eqref{eq:Acl}, and $u\in\mbb{U}$, $y\in\mbb{Y}$.
  \end{enumerate} for all $k\in\Nset$
  \label{prop:dd-oinf}
\end{proposition}
\begin{proof}
  $i)$ $\Xi_0$ is compact with the origin in its relative interior (Proposition~\ref{prop:ZAAR}) and $\mbb{V}$ is $A_{\mn{cl}}$-invariant (Lemma~\ref{lem:history}); finite determination follows from \cite[Thm.~4.1]{Gilbert1991} applied on $\mbb{V}$. $ii)$ follows from Proposition~\ref{prop:ZAAR}.
\end{proof}

\section{Feasible sets and invariance}
\label{sec:feas_sets}
In this section, we study and analyse the feasible sets of the \ac{DeePC} problem and their invariance properties. By Lemma~\ref{lem:QP_cond}, the histories for which \eqref{eq:DeePC_QP} is feasible form the \emph{feasible set}
\begin{equation}
  \mc{F}_N := \{\xi\colon \exists\,\mb{u}\in\Rset^{mN},~G\mb{u}\leq S_\xi\xi+w\},
  \label{eq:FN}
\end{equation}
the projection onto $\xi$ of a polyhedron in $\Rset^{(m+p)T_\ini+mN}$. The \ac{DeePC} feasible set \eqref{eq:FN} depends on the horizon length $N$ and is not necessarily an invariant set for small $N$ because its underlying \ac{DeePC} problem lacks terminal ingredients \cite{Berberich2021a}. As $N$ grows, the feasible set approaches the maximal control invariant, as shown in the model-based case in \cite{SchulzeDarup2014, Liu2022}.

\subsection{Critical Regions and Affine Control Laws}
\label{sec:dd-pwa-regions}

When the activation pattern of \eqref{eq:NN_implicit} is fixed, \eqref{eq:NN_law} is affine in $\xi$; the following result identifies the set on which each pattern is valid.

\begin{proposition}[Critical regions]
  Suppose Assumptions~\ref{assump:Tini} and \ref{assump:constraints} hold. Fix $\mb{r}=0$. Let $\mf{A} = \{\mc{A}\subseteq\{1,\ldots,n_c\}\colon G_\mc{A}\ts{ has full row rank}\}$, for $\mc{A}\in\mf{A}$, let $\lambda_\mc{A}(\xi) := -M_{\mc{A}\mc{A}}^{-1}(\hat S_\mc{A}\xi+w_\mc{A})$ and
  \begin{equation}
    \mc{R}_\mc{A} := \{\xi\colon \lambda_\mc{A}(\xi)\geq0,~\hat S_{\bar{\mc{A}}}\xi+w_{\bar{\mc{A}}}+M_{\bar{\mc{A}}\mc{A}}\lambda_\mc{A}(\xi)\geq0\}.
    \label{eq:crit-region}
  \end{equation}
  \begin{enumerate}[i)]
  \item For $\xi\in\mc{R}_\mc{A}$, the vector $\zeta$ with $\zeta_\mc{A} = \lambda_\mc{A}(\xi)$ and $\zeta_{\bar{\mc{A}}} = -(\hat S_{\bar{\mc{A}}}\xi+w_{\bar{\mc{A}}}+M_{\bar{\mc{A}}\mc{A}}\lambda_\mc{A}(\xi))$ solves \eqref{eq:NN_implicit}.
    \item On $\mc{R}_\mc{A}$, $\mb{u}^*(\xi) = K_\mc{A}\xi+c_\mc{A}$ with $K_\mc{A} := -H^{-1}(F_\xi^\top-G_\mc{A}^\top M_{\mc{A}\mc{A}}^{-1}\hat S_\mc{A})$ and $c_\mc{A} := H^{-1}G_\mc{A}^\top M_{\mc{A}\mc{A}}^{-1}w_\mc{A}$.
    \item $\mc{F}_N = \bigcup_{\mc{A}\in\mf{A}}\mc{R}_\mc{A}$.
  \end{enumerate}
  \label{prop:crit}
\end{proposition}
\begin{proof}
  $i)$ Substitute $\phi(\zeta)$, with $\phi(\zeta)_\mc{A} = \lambda_\mc{A}(\xi)\geq0$ and $\phi(\zeta)_{\bar{\mc{A}}} = 0$, in \eqref{eq:NN_implicit} and use $M_{\mc{A}\mc{A}}\lambda_\mc{A}(\xi) = -(\hat S_\mc{A}\xi+w_\mc{A})$. $ii)$ Follows from \eqref{eq:NN_law} with $W_f\phi(\zeta) = -H^{-1}G_\mc{A}^\top\lambda_\mc{A}(\xi)$. $iii)$ $\supseteq$ follows from $i)$ and Lemma~\ref{lem:mpc-to-nn}. The converse $\subseteq$, for $\xi\in\mc{F}_N$, $-(H\mb{u}^*+F_\xi^\top\xi)$ is a nonnegative combination of the rows of $G$ active at $\mb{u}^*$; by Carath\'eodory's theorem it is a nonnegative combination of linearly independent active rows indexed by some $\mc{A}\in\mf{A}$, and then $\xi\in\mc{R}_\mc{A}$.
\end{proof}

Proposition~\ref{prop:crit} is the data-driven counterpart of the explicit solution of multiparametric QPs \cite{Bemporad2002} because we obtain a similar partition of the feasible set $\mc{F}_N$ into linear regions. This proposition includes the results in Proposition~\ref{prop:ZAAR}, \ie $\mc{R}_{\emptyset} = \Xi_0$ and $K_\emptyset = -H^{-1}F_\xi$. By Lemma~\ref{lem:history}, on $\mc{R}_\mc{A}$ the closed loop is affine,
\begin{equation}
  \xi^+ = (A_\xi + B_\xi E_1K^{\mb{u}}_\mc{A})\xi + B_\xi E_1c^{\mb{u}}_\mc{A}.
  \label{eq:pwa-history-map}
\end{equation}
We denote by $f_\pi$ the resulting continuous \ac{PWA} map on $\mc{F}_N$.
\subsection{Region of recursive feasibility}
\label{sec:dd-pwa-invariance}
Feasibility of the QP \eqref{eq:DeePC_QP} at $\xi$ does not imply feasibility at $\xi^+$, as no terminal ingredient is imposed \cite{Berberich2021a}. Using the \ac{MPI} algorithm with initialisation $\mc{S}_0 := \mc{X}\cap\mc{F}_N$, the \ac{DeePC} feasible set is
\begin{equation}
  \mc{C}_\pi := \mc{O}(f_\pi,\mc{S}_0).
  \label{eq:Cpi-dd}
\end{equation}
For $\mc{S}\subseteq\mc{F}_N$, define the predecessor under the \ac{DeePC} law
\begin{equation}
  \Pre_\pi(\mc{S}) := \textstyle\bigcup_{\mc{A}}\{\xi\in\mc{R}_\mc{A}\colon A_\mc{A}\xi+a_\mc{A}\in\mc{S}\}
  \label{eq:pwa-pre}
\end{equation}
and the nested sequence $\mc{S}_{j+1} := \mc{S}_j\cap\Pre_\pi(\mc{S}_j)$.
\begin{proposition}[Recursive feasibility]
  Suppose Assumptions~\ref{assump:Tini}--\ref{assump:constraints} hold and $\mb{r}=0$.
  \begin{enumerate}[i)]
  \item $\xi\in\mc{C}_\pi$ iff $\xi\in\mc{X}$ and QP \eqref{eq:DeePC_QP} is feasible along the \ac{DeePC} closed loop trajectory.
  \item $\forall j$, $\mc{S}_j\supseteq\mc{S}_{j+1}\supseteq\mc{C}_\pi$.
  \end{enumerate}
  \label{prop:Cpi}
\end{proposition}
\begin{proof}
  $i)$ If QP \eqref{eq:DeePC_QP} is feasible at $\xi$, the applied input and, by Lemma~\ref{lem:history}, the measured output coincide with the first predicted pair, which satisfies \eqref{eq:DeePC_u}--\eqref{eq:DeePC_y}; hence $\xi\in\mc{X}$ implies $\xi^+\in\mc{X}$ as in Proposition~\ref{prop:dd-oinf}. $ii)$ $\mc{C}_\pi\subseteq\Pre_\pi(\mc{C}_\pi)$ and induction give $\mc{C}_\pi\subseteq\mc{S}_j$. If $\mc{S}_{j+1} = \mc{S}_j$, then $\mc{S}_j\subseteq\Pre_\pi(\mc{S}_j)$, so $\mc{S}_j\subseteq\mc{S}_0$ is positively invariant and $\mc{S}_j\subseteq\mc{C}_\pi$.
\end{proof}


\subsection{Invariant set hierarchy}
\label{sec:invar-set-hier}
By Lemma~\ref{lem:history}, the \ac{MCI} subset of $\mc{X}$,
\begin{equation}
  \mc{C}_\infty := \{\xi\in\mc{X}\colon \exists \{u_k\}_{k\in\Nset}\subset\mbb{U},~A_\xi\xi+B_\xi u_k\in\mc{X}~\forall k\in\Nset\},
  \label{eq:Cinf-dd}
\end{equation}
is computed from data by $\Omega_0 := \mc{X}$, $\Omega_{j+1} := \Omega_j\cap\Pre_c(\Omega_j)$, where $\Pre_c(\Omega) := \{\xi\colon\exists u\in\mbb{U},~A_\xi\xi+B_\xi u\in\Omega\}$ is a projection that eliminates $m$ variables \cite{Bertsekas1972,Blanchini1999,Blanchini2008}; \cite{Bajelani2025} constructs such sets in input/output coordinates. The following is our main result and clarifies the relationship between the sets computed in the previous sections.
\begin{theorem}[Set hierarchy]
  Suppose Assumptions~\ref{assump:Tini}--\ref{assump:constraints} hold and $\mb{r} = 0$. Then
  \begin{equation} \mc{O}_\infty\subseteq\mc{C}_\pi\subseteq\mc{C}_\infty\subseteq\mc{S}_0\subseteq\mc{F}_N,\qquad \mc{O}_\infty\subseteq\Xi_0\subseteq\mc{S}_0.
    \label{eq:set-hierarchy}
  \end{equation}
  \label{thm:hierarchy}
\end{theorem}
\begin{proof}
  On $\Xi_0$, $f_\pi = A_{\mathrm{cl}}$; hence $\mc{O}_\infty$ is positively invariant under $f_\pi$ and $\mc{O}_\infty\subseteq\Xi_0\subseteq\mc{S}_0$, so $\mc{O}_\infty\subseteq\mc{C}_\pi$ by maximality. By Proposition~\ref{prop:Cpi}~i), $\mc{C}_\pi\subseteq\mc{X}$ is invariant under $\xi^+ = A_\xi\xi+B_\xi u$ with $u = E_1\mb{u}^*(\xi)\in\mbb{U}$, so $\mc{C}_\pi\subseteq\mc{C}_\infty$. Let $\xi_0\in\mc{C}_\infty$ with admissible inputs $\{u_k\}_{k}$; the outputs $y_0,\ldots,y_{N-1}$ are entries of $\xi_1,\ldots,\xi_N\in\mc{X}$, so $(u_0,\ldots,u_{N-1})$ is feasible for \eqref{eq:DeePC_QP} by Lemma~\ref{lem:QP_cond}, \ie $\xi_0\in\mc{F}_N$. The remaining inclusions hold by definition.
\end{proof}

The inclusion chains in \eqref{eq:set-hierarchy} are strict, \ie $\mc{S}_0\setminus\mc{C}_\infty$ contains histories at which \ac{DeePC} is feasible although no admissible input avoids a future violation. The horizon $N$ plays an important role, as $\mc{F}_{N+1}\subseteq\mc{F}_N$. Similarly, $\mc{C}_\infty\setminus\mc{C}_\pi$ shrinks as $N$ grows. Our last result addresses how positive invariant sets $\mc{P}$ from a linear gain $K$, \ie certified sets by the linear programs of \cite{Bisoffi2020, Bisoffi2023} applied to \eqref{eq:history-model}, belong to the hierarchy of Theorem~\ref{thm:hierarchy}.

\begin{corollary}[Linear-feedback certificates]
  \label{cor:linear}
  Suppose Assumptions~\ref{assump:Tini}--\ref{assump:constraints} hold and $\mb{r}=0$. For $K\in\Rset^{m\times(m+p)T_\ini}$, let $\mc{M}(K) := \mc{O}\big(A_\xi+B_\xi K,\ \mc{X}\cap\{\xi\colon K\xi\in\mbb{U}\}\big)$. Then $\mc{O}_\infty\subseteq\mc{M}(K_\emptyset)$ and $\mc{M}(K)\subseteq\mc{C}_\infty$ for every $K$. In particular, every $\mc{P}\subseteq\mc{X}$ that is positively invariant for $\xi^+ = (A_\xi+B_\xi K)\xi$ with $K\xi\in\mbb{U}$ for all $\xi\in\mc{P}$ satisfies $\mc{P}\subseteq\mc{M}(K)\subseteq\mc{C}_\infty$.
\end{corollary}
\begin{proof}
  $\mc{M}(K)\subseteq\mc{X}$ is positively invariant under the admissible input $u = K\xi$, hence $\mc{M}(K)\subseteq\mc{C}_\infty$ by \eqref{eq:Cinf-dd}. The set $\mc{O}_\infty\subseteq\Xi_0\subseteq\mc{X}$ is positively invariant under $A_{\mathrm{cl}} = A_\xi+B_\xi K_\emptyset$, and $K_\emptyset\xi = E_1\mb{u}^*(\xi)\in\mbb{U}$ on $\Xi_0\subseteq\mc{R}_\emptyset$ by \eqref{eq:cond_cons}; maximality of $\mc{M}(K_\emptyset)$ gives $\mc{O}_\infty\subseteq\mc{M}(K_\emptyset)$, and maximality of $\mc{M}(K)$ gives $\mc{P}\subseteq\mc{M}(K)$.
\end{proof}

Corollary~\ref{cor:linear} adds a branch to \eqref{eq:set-hierarchy}: $\mc{O}_\infty\subseteq\mc{M}(K_\emptyset)\subseteq\mc{C}_\infty$ for the gain of the unsaturated branch and $\mc{M}(K)\subseteq\mc{C}_\infty$ for every other gain, while no inclusion between $\mc{M}(K)$ and $\mc{C}_\pi$ follows from Theorem~\ref{thm:hierarchy}. Hence $\mc{P}\subseteq\mc{C}_\infty$ is necessary for the linear programs of \cite{Bisoffi2020,Bisoffi2023} to be feasible for a candidate set $\mc{P}$, whatever the gain, and $\mc{C}_\infty$ bounds every certifiable set. The converse fails: a polyhedral control invariant set need not be invariant under any linear feedback, as noted in \cite[Rem.~4]{Bisoffi2020} and by the definition of $\mc{C}_\infty$. The programs of \cite{Bisoffi2020,Bisoffi2023} certify a given set, robustly to noise in \cite{Bisoffi2023}, with one linear program whose size grows with $T$; the recursions above compute maximal sets for exact data, at a cost that does not depend on $T$ once $g$ is eliminated but is exponential in $n_c$ and $mN$ in the worst case.

\section{Illustrative examples}
\label{sec:illustr-exampl}
\begin{figure}[!t]
  \centering
  \captionsetup[subfloat]{
    position=top,
    justification=raggedright,
    singlelinecheck=false,
    captionskip=-15pt 
  }
  \subfloat[]{
%
%
\begin{tikzpicture}

\begin{axis}[%
width=0.4\textwidth,
height=0.15\textwidth,
at={(0,0)},
scale only axis,
xmin=-1.2,
xmax=1.2,
xlabel={\footnotesize $x_1$},
xmajorgrids,
ymin=-2,
ymax=2,
ylabel={\footnotesize $x_2$},
ymajorgrids,
legend style={at={(0.03,0.03)},
  anchor=south west,
  legend cell align=left,
  align=left,
  draw=none,
  fill=none,
  font=\footnotesize
}
]
        %
        %
\addplot [area legend,color=i_col_1,fill = i_col_1!50,fill opacity = 0.5,solid,line width=1pt]
  table[row sep=crcr]{%
-1.0000000016	0.4612281674\\
-1	0.3377192981\\
-0.5826881305	-0.2790969916\\
0.4188612648	-1.7594572396\\
0.4188612759	-1.7594572369\\
0.7875000015	-1.2145833331\\
0.9999999999	-0.5416666672\\
1	-0.4386043176\\
0.9999999999	-0.337719298\\
-0.4188612759	1.7594572368\\
-0.418861276	1.7594572368\\
-0.7874999981	1.2145833389\\
-0.9999999976	0.5416666736\\
-1.0000000016	0.4612281674\\
}--cycle;
\addlegendentry{$\Xi_0$};


\addplot[area legend,color=col_1,fill opacity = 0.2,fill = col_1!50!white,solid,line width=1.0pt,dashed]
  table[row sep=crcr]{%
-1.0000000001	0.4266693409\\
-1	0.3377192982\\
0.0600438522	-1.2290999045\\
0.5395097396	-1.5083013349\\
0.608978574	-1.478450521\\
0.6089785814	-1.4784505178\\
0.6988325003	-1.3456401208\\
0.7875	-1.2145833333\\
1	-0.5416666667\\
1	-0.3983291593\\
0.9999999983	-0.3377192935\\
-0.0600438523	1.2290999044\\
-0.5395097397	1.5083013349\\
-0.6089785814	1.4784505178\\
-0.7875	1.2145833333\\
-1	0.5416666666\\
-1.0000000001	0.4266693409\\
}--cycle;
\addlegendentry{$\mathcal{O}_\infty$};

\end{axis}
\end{tikzpicture}%
}\\\vspace{-0.2cm}
  \subfloat[]{
%
\definecolor{mycolor1}{rgb}{0.20000,0.40000,0.80000}%
\definecolor{mycolor2}{rgb}{0.10000,0.60000,0.30000}%
\definecolor{mycolor3}{rgb}{0.80000,0.15000,0.15000}%
\begin{tikzpicture}

\begin{axis}[%
width=0.4\textwidth,
height=0.15\textwidth,
at={(0\textwidth,0\textwidth)},
scale only axis,
xmin=-1.2,
xmax=1.2,
xlabel={$x_1$},
xmajorgrids,
ymin=-2,
ymax=2,
ylabel={$x_2$},
ymajorgrids,
legend style={at={(0.03,0.03)},
  anchor=south west,
  legend cell align=left,
  align=left,
  draw=none,
  fill=none,
  font=\footnotesize
}
]
\addplot [area legend,fill = i_col_10!50,fill opacity = 0.2,color=i_col_10,solid,line width=1pt]
  table[row sep=crcr]{%
-1	0.3377192982\\
0.0600438524	-1.229099904\\
0.493678078	-1.4816127511\\
0.5395097397	-1.5083013349\\
0.7145960422	-1.4330666098\\
0.9839080493	-1.3173430222\\
1	-1.264507925\\
1	-0.7388955245\\
1	-0.3377192982\\
-0.0600438521	1.2290999048\\
-0.288610647	1.3621983843\\
-0.5395097379	1.5083013366\\
-0.5395097429	1.5083013354\\
-0.9839080456	1.3173430303\\
-0.9993568465	1.2666196104\\
-1	1.264507925\\
-1	0.7388955251\\
-1	0.3377192982\\
}--cycle;
\addlegendentry{$\mathcal{O}_\infty^{model}$};


\addplot [area legend,color=col_1,fill = col_1!50,fill opacity = 0.2,solid,line width=1pt,dashed]
  table[row sep=crcr]{%
-1.0000000001	0.4266693409\\
-1	0.3377192982\\
0.0600438522	-1.2290999045\\
0.5395097396	-1.5083013349\\
0.608978574	-1.478450521\\
0.6089785814	-1.4784505178\\
0.6988325003	-1.3456401208\\
0.7875	-1.2145833333\\
1	-0.5416666667\\
1	-0.3983291593\\
0.9999999983	-0.3377192935\\
-0.0600438523	1.2290999044\\
-0.5395097397	1.5083013349\\
-0.6089785814	1.4784505178\\
-0.7875	1.2145833333\\
-1	0.5416666666\\
-1.0000000001	0.4266693409\\
}--cycle;
\addlegendentry{$\mathcal{O}_\infty$};

\end{axis}
\end{tikzpicture}%
}\\\vspace{-0.2cm}
  \subfloat[]{
%
\definecolor{mycolor1}{rgb}{0.20000,0.40000,0.80000}%
\definecolor{mycolor2}{rgb}{0.10000,0.60000,0.30000}%
\definecolor{mycolor3}{rgb}{0.80000,0.15000,0.15000}%
\begin{tikzpicture}

\begin{axis}[%
width=0.4\textwidth,
height=0.15\textwidth,
at={(0\textwidth,0\textwidth)},
scale only axis,
xmin=-1.7,
xmax=1.2,
xlabel={$x_1$},
xmajorgrids,
ymin=-4,
ymax=4,
ylabel={$x_2$},
ymajorgrids,
view={0}{90},
legend style={at={(0.03,0.03)},
  anchor=south west,
  legend cell align=left,
  align=left,
  draw=none,
  fill=none,
  font=\footnotesize}
]

\addplot3[area legend,dotted,
line  width=1.0pt,
fill=col_2!50,fill opacity=0.2,
draw=col_2]
table[row sep=crcr] {%
x	y	z\\
1	-3.4956140351	0\\
-0.7792207792	-0.8658008658	0\\
-1	-0.1666666667	0\\
-1	3.4956140351	0\\
0.7792207792	0.8658008658	0\\
1	0.1666666667	0\\
}--cycle;
\addlegendentry{$\mc{F}_N$};

\addplot3[area legend,solid,
line  width=1.0pt,
fill=col_3!50,fill opacity=0.2,
draw=col_3]
table[row sep=crcr] {%
x	y	z\\
1.0000000182	-0.5416666243	0\\
0.7875000089	-1.2145833144	0\\
0.7874999935	-1.2145833482	0\\
0.319704264	-1.9060182697	0\\
0.2037202877	-2.0774507263	0\\
0.122125371	-2.198053728	0\\
0.1221253505	-2.1980537247	0\\
-0.5034088502	-1.2734703418	0\\
-0.7792207798	-0.8658008647	0\\
-0.9999999929	-0.1666666915	0\\
-1	-0.0683510214	0\\
-1.0000000005	0.5416666429	0\\
-0.8025820684	1.1668234515	0\\
-0.7875000009	1.2145833298	0\\
-0.7874999976	1.2145833368	0\\
-0.3243069606	1.8992152007	0\\
-0.1221253836	2.1980537305	0\\
-0.1221253709	2.1980537281	0\\
0.7792207792	0.8658008658	0\\
1	0.1666666669	0\\
}--cycle;
\addlegendentry{$\mc{S}_0$};

\addplot3[area legend,solid,
line  width=1.0pt,
fill=black!50,fill opacity=0.2,
draw=black]
table[row sep=crcr] {%
x	y	z\\
1.0000000001	-0.3377192986	0\\
0.9999999999	-0.5416666672	0\\
0.7875000015	-1.2145833331	0\\
0.2441360357	-2.0177134033	0\\
0.244136035	-2.017713404	0\\
-0.3433968629	-1.5099791983	0\\
-0.7792197585	-0.8658023746	0\\
-0.7792207792	-0.8658008658	0\\
-0.936431989	-0.367965369	0\\
-0.9999999996	-0.1666666676	0\\
-1	-0.0815404505	0\\
-1.0000000016	0.4612281674	0\\
-0.9999999976	0.5416666736	0\\
-0.7874999981	1.2145833389	0\\
-0.2441360353	2.0177134037	0\\
0.3433968629	1.5099791983	0\\
0.7792207811	0.8658008584	0\\
0.9364319891	0.367965368	0\\
0.9999999993	0.1666666682	0\\
1	0.0801310516	0\\
}--cycle;
\addlegendentry{$\mc{C}_\pi$};

\addplot3[area legend,dash pattern=on 1pt off 3pt on 3pt off 3pt,
line  width=1.0pt,
fill=col_4!50,fill opacity=0.2,
draw=col_4]
table[row sep=crcr] {%
x	y	z\\
1	-0.5416666667	0\\
0.7875000001	-1.2145833332	0\\
0.2441360545	-2.017713381	0\\
0.2441360353	-2.017713404	0\\
0.2441360313	-2.0177134025	0\\
-0.3433968629	-1.5099791983	0\\
-0.7792207792	-0.8658008658	0\\
-0.9999999998	-0.1666666673	0\\
-1	0.1211915225	0\\
-1	0.5416666667	0\\
-0.7875000001	1.2145833331	0\\
-0.7875	1.2145833333	0\\
-0.4672249633	1.6879723717	0\\
-0.2441359893	2.0177135063	0\\
0.3433967196	1.5099793221	0\\
0.3433968628	1.5099791983	0\\
0.7792207792	0.8658008658	0\\
1	0.1666666667	0\\
}--cycle;
\addlegendentry{$\mc{C}_\infty$};

\addplot3[area legend,dashed,
line  width=1.0pt,
fill=col_1!50,fill opacity=0.2,
draw=col_1]
table[row sep=crcr] {%
x	y	z\\
0.9999999983	-0.3377192935	0\\
1	-0.3983291593	0\\
1	-0.5416666667	0\\
0.7875	-1.2145833333	0\\
0.6988325003	-1.3456401208	0\\
0.6089785814	-1.4784505178	0\\
0.608978574	-1.478450521	0\\
0.5395097396	-1.5083013349	0\\
0.0600438522	-1.2290999045	0\\
-1	0.3377192982	0\\
-1.0000000001	0.4266693409	0\\
-1	0.5416666666	0\\
-0.7875	1.2145833333	0\\
-0.6089785814	1.4784505178	0\\
-0.5395097397	1.5083013349	0\\
-0.0600438523	1.2290999044	0\\
}--cycle;
\addlegendentry{$\mc{O}_\infty$};

\end{axis}
\end{tikzpicture}%
}
  \caption{(a) $\Xi_0$, the iterate $\mc{P}_1$ of \eqref{eq:gt-dd} and $\mc{O}_\infty = \mc{P}_2$. (b) $\mc{O}_\infty$, $\mc{O}(A_{\mathrm{cl}},\mc{R}_\emptyset)$ (dashed) and the invariant set of the model-based \ac{MPC}' s unsaturated branch.}
  \label{fig:inv_deepc}
\end{figure}
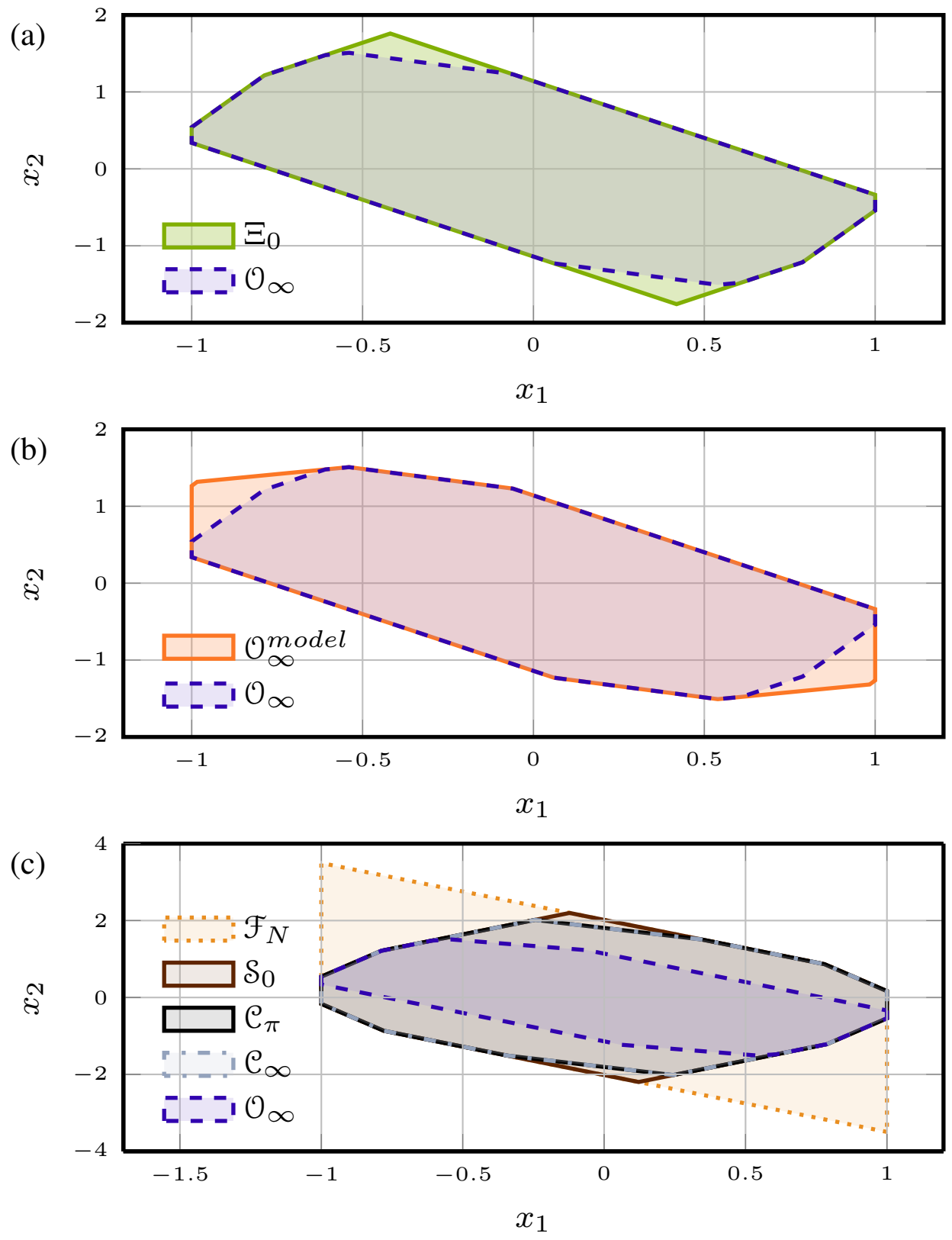

Consider the following baseline system
\[ x^+=\begin{bmatrix} 0.95 & 0.3\\ -0.2 & 0.95 \end{bmatrix} x+ \begin{bmatrix} 0\\ 1 \end{bmatrix}u, \qquad y= \begin{bmatrix} 1 & 0 \end{bmatrix} x, \]
with $n = \ell = 2$. A noise-free trajectory of length $T = 500$ was generated from $x_0 = 0$. With $T_\ini = 2$ and $N=3$, $L = 5$ and $\rank\mc{H}_{L+n}(u^d) = 7$, so Assumption~\ref{assump:Tini}  holds. The constraints are $|u|\leq 0.5$ and $|y|\leq 1$, and $Q = 1$, $R = 0.3$, so $n_c = 12$. Since $pT_\ini = n$, $\mbb{V} = \Rset^4$. Every set in Figure~\ref{fig:inv_deepc} is computed from the Hankel matrices only. For display, histories are mapped to the current state $x_k$ through the true model; as $\xi\mapsto(x,u_\ini)$ is a linear bijection, each figure shows $\pi(\mc{P}) := \{x\colon\xi\in\mc{P}\}$.

In Figure~\ref{fig:inv_deepc} (a), we observe the relation between the unsaturated branch $\Xi_0\cap\mc{X}$ and the largest invariant set contained within it $\mc{O}_\infty$. For this system, the algorithm in \cite{Gilbert1991} finished after two iterations, which explains the similar size of these two sets. The data-driven set implicitly constrains inputs and outputs, which makes for more stringent constraints which yield the difference seen in Figure \ref{fig:inv_deepc} (b). The remaining difference is only because of $\mc{X}$: $\pi(\mc{O}_\infty) = \pi(\mc{X})\cap\pi(\mc{O}(A_{\mathrm{cl}},\mc{R}_\emptyset))$, where $\pi(\mc{X})$ is the set of states reachable through an admissible window. The two removed pieces, of area $0.118$ each, contain states such as $x_k = (-0.984,\,1.317)$, which can be reached with $|y_{k-1}|\leq 1$ only if $u_{k-1}\geq 1.225$. Lastly, for this example, the \ac{PWA} affine control law yielded three regions; the recursion of Proposition~\ref{prop:Cpi}~$iii)$ terminated after two iterations, certifying $\mc{C}_\pi$, see Figure~\ref{fig:inv_deepc} (c). In $\Rset^4$, $\mc{O}_\infty\subseteq\mc{C}_\pi\subseteq\mc{C}_\infty\subseteq\mc{S}_0$ and $\Xi_0\subseteq\mc{C}_\infty$ were verified, and $\mc{C}_\infty\setminus\mc{C}_\pi$ has empty interior; since both sets are closed, $\mc{C}_\pi = \mc{C}_\infty$. Thus, for this design, the \ac{DeePC} policy loses nothing with respect to an unrestricted admissible input, and its region of recursive feasibility is $97.02\%$ of $\mc{S}_0$ in projection; the remaining $2.98\%$ is the cost of the horizon $N=3$ without terminal ingredients. Restricting \ac{DeePC} to its unsaturated branch retains $60.72\%$ of $\mc{C}_\pi$.

\section{Conclusions and future work}
\label{sec:concl-future-work}

The \ac{DeePC} quadratic program admits an exact condensed form on the consistent histories, whose solution is an implicit \ac{NN} with weights built from Hankel matrices. Its activation patterns are the critical regions of the QP, which yield the feasible set of \ac{DeePC}, an exact linear model of the history window, and three nested invariant sets, all computed from the Hankel matrices. In the example, the region of recursive feasibility of \ac{DeePC} without terminal ingredients coincides with the maximal control invariant set of admissible histories. The construction requires exact data, and enumerating activation patterns grows combinatorially with $n_c$ and $mN$. Future work will address noisy data, where the condensed problem becomes a least-squares predictor, the sets must be robustified, and $\mc{C}_\pi$ and $\mc{C}_\infty$ can be used as data-driven terminal ingredients.

\bibliography{biblio}

\end{document}